\documentclass{article}
\usepackage{amsfonts,amssymb,latexsym,amsmath,amscd,euscript}
\newcounter{num}[section]
\newcommand{\Num}{\refstepcounter{num}%
	\textbf{\arabic{section}.\arabic{num}}}

\newcommand{\Theorem}{\textbf{Theorem~}}
\newcommand{\Proof}{{\LARGE\emph{Proof}}}
\newcommand{\Def}{\textbf{Definition}}

\newcommand{\Prop}{\textbf{Proposition~}}
\newcommand{\Cor}{ \textbf{ Corollary~}}

\newcommand{\GL}{{\mathrm{GL}}}
\newcommand{\SL}{{\mathrm{SL}}}
\renewcommand{\leq}{\leqslant}

\newcommand{\Ad}{{\mathrm{Ad}}}
\newcommand{\UT}{{\mathrm{UT}}}

\newcommand{\Pb}{{\Bbb P}}
\newcommand{\Db}{{\Bbb D}}
\newcommand{\Bb}{{\Bbb B}}
\newcommand{\Nb}{{\Bbb N}}
\newcommand{\Yb}{{\Bbb Y}}

\newcommand{\al}{\alpha}

\newcommand{\Ax}{{\mathfrak A}}
\newcommand{\Mat}{{\mathrm{Mat}}}

\newcommand{\Hc}{{\mathcal H}}

\newcommand{\Fc}{{\mathcal F}}
\newcommand{\Sc}{{\mathcal S}}

\newcommand{\Pc}{{\mathcal P}}
\newcommand{\Qc}{{\mathcal Q}}

\begin{document}
\Large

\title{Fields of $U$-invariants of matrix tuples }
\author{A.N. Panov\footnote{The work is supported by the RFBR grant   20-01-00091a}}
\date{}
 \maketitle
 {\small   Mechanical and Mathematical  Department, Samara National Research University, Samara, Russia} 
 \begin{abstract}
 The general linear group $\GL(n)$ acts on the direct sum of $m$ copies of  $\Mat(n)$ by the adjoint action. The action of $\GL(n)$ induces the action of the unitriangular subgroup $U$. 
  We  present the system of free generators of the field of $U$-invariants. 
 	 \end{abstract}
{\small Keywords: theory of invariants, adjoint representation, matrix tuple}\\
 {\small 2010 Mathematics Subject Classifications: 15A15, 15A72, 13A50}

\section{Introduction}  
One of the classical problems of the theory of invariants is the problem of description of invariants of matrix tuples. Let $K$ be an arbitrary field. Let  $\Mat(n)$ stand for the linear space of  $(n\times n)$-matrices with entries in the field  $K $ and  $\GL(n)$ be the general linear group of order  $n$.
Let us  consider the linear space   
$$\Hc=\Mat(n)\oplus\ldots\oplus \Mat(n)$$ of   $m$-tuples   of  $\Mat(n)$. The group  $\GL(n)$ acts on  $\Hc$ by the formula 
$$ \Ad_g(X_1,\ldots, X_m) = (gX_1 g^{-1},\ldots, g X_mg^{-1}),$$
where
 $g\in \GL(n)$ 
and $X_1,\ldots, X_m\in \Mat(n)$. 
The action of the group  $\GL(n)$ on $\Hc$ defines the representation  $$\rho(g)f(X_1,\ldots,X_m) = f(g^{-1}X_1g,\ldots, g^{-1}X_mg)$$ of the group  $\GL(n)$ in the space of regular functions $K[\Hc]$. This representation is extended to the action of   $\GL(n)$ on the field of rational functions  $K(\Hc)$.
For a given subgroup $G\subseteq GL(n)$, the problem is to describe the algebra (respectively, the field) of invariants with respect to the action of $G$ on $\Hc$. 

In the case  $G=\GL(n)$ (or $G=\SL(n)$),  this problem is solved in the framework of the classical theory of invariants in tensors  (see \cite{PV, Pro, Don}). The algebra of $\GL(n)$-invariants is generated by the system of polynomials  $\mathrm{Tr}(A_{i_1}\cdots A_{i_p})$, where $1\leq i_1,\ldots,i_p\leq m$.

The group  $\GL(n)$ contains the subgroup of unitriangular matrices  $U=\UT(n)$, which consists of the upper triangular matrices with ones on the diagonal.  We consider the algebra of  $U$-invariants  $K[\Hc]^U$ and the field of  $U$-invariants $K(\Hc)^U$. 

Little is known about the structure of the algebra of invariants   $K[\Hc]^U $.  It follows from \cite[Theorem 3.13]{PV} that the algebra  $K[\Hc]^U $ is finitely generated.   For  $K[\Hc]^U $, the problem of construction of system of generators with their relations  is an unsoled problem even for  $m=1$. 

The structure of  the field $K(\Hc)^U$ is much simpler. The  group  $U$ acts on  $\Hc$ by unipotent transformations. It implies that\\ 
1) the field   $K(\Hc)^U$ is rational \cite{Mi},  i.e. it is a  pure transcendental extension over the main field  $K$,\\
2) the field $K(\Hc)^U$ is a field of fractions of the algebra of invariants  $K[\Hc]^U$ ~\cite[Theorem 3.3]{PV}.

The goal of this paper is to construct a system of free generators of the field of  $U$-invariants $K(\Hc)^U$.  
For the adjoint action of  $U$ on $\Mat(n)$ (the case $m=1$), the system generators of  the field of $U$-invariants is presented in Theorem   \ref{ThUone}. 
This system of generators is not unique. The another system of free  generators of $K(\Mat(n))^U$ was constructed in  \cite{PV1, Pan}.  

For an arbitrary $m$, the system of free generators of $K(\Hc)^U$ is presented in Theorem \ref{ThU}.   
We apply invariants  $\{P_{ik}(X,Y)\}$ from (\ref{PXY}).  The similar invariants  were earlier used in the papers  \cite{PSev, Sev}.

\section{Field of  $U$-invariants on  $\Mat(n)$}

For $m=2$, we have $\Hc_2=\Mat(n)\oplus \Mat(n)$.
Let $\{x_{ij}\}_{i,j=1}^n$ and  $\{y_{ij}\}_{i,j=1}^n$  be the systems of standard coordinate  functions on the first and second components of $\Hc_2$.  
Consider two matrices  
$$X=(x_{ij})_{i,j=1}^n\qquad\mbox{and}\qquad Y=(y_{ij})_{i,j=1}^n.$$ 

For two positive integers  $a$ and $b$,  we denote by $[a,b]$ the subset of integers  $a\leq i\leq b$. 
For any integer   $1\leq i\leq n$,  let  $i'$ be the symmetric number to  $i$ with respect to the  center of the segment  $[1,n]$. We have $i'=n-i+1$.

For the pair  $i'\leq j$ (i.e.  $(i,j)$ lies on or below the anti-diagonal), let  $M_{ij}(X)$ be the minor of order  $i'$ of the matrix   $X$ with the system of rows $[i,n]$ and columns  $[1,i'-1]\sqcup \{j\}$.

For the pair $j\leq k$, let  $N_{jk}(Y)$  be the minor of order  $k'$ of the matrix   $Y$ with the system of rows  $\{j\}\sqcup [k+1,n]$ and columns $[1,k']$.\\
\\
{\bf Example}. For $n=5$, we have
{\small $$X=\left(\begin{array}{ccccc}
x_{11}&x_{12}&x_{13}&x_{14}&x_{15}\\
x_{21}&x_{22}&x_{23}&x_{24}&x_{25}\\
x_{31}&x_{32}&x_{33}&x_{34}&x_{35}\\
x_{41}&x_{42}&x_{43}&x_{44}&x_{45}\\
x_{51}&x_{52}&x_{53}&x_{54}&x_{55}\\
\end{array}\right),\quad M_{34}(X) = \left|\begin{array}{ccc}
x_{31}&x_{32}&x_{34}\\
x_{41}&x_{42}&x_{44}\\
x_{51}&x_{52}&x_{54}
\end{array}\right|, \quad M_{45}(X)=\left|\begin{array}{cc}x_{41}&x_{45}\\
x_{51}&x_{55}\end{array}\right|,\quad M_{53}(X)=x_{53},$$}
{\small$$Y=\left(\begin{array}{ccccc}
	y_{11}&y_{12}&y_{13}&y_{14}&y_{15}\\
	y_{21}&y_{22}&y_{23}&y_{24}&y_{25}\\
	y_{31}&y_{32}&y_{33}&y_{34}&y_{35}\\
	y_{41}&y_{42}&y_{43}&y_{44}&y_{45}\\
	y_{51}&y_{52}&y_{53}&y_{54}&y_{55}\\
	\end{array}\right),\quad N_{23}(Y) = \left|\begin{array}{ccc}
	y_{21}&y_{22}&y_{23}\\
	y_{41}&y_{42}&y_{43}\\
	y_{51}&y_{52}&y_{53}
	\end{array}\right|, \quad N_{14}(Y)=\left|\begin{array}{cc}y_{11}&y_{12}\\
	y_{51}&y_{52}\end{array}\right|, \quad N_{35}(Y) = y_{31}.$$}

Let $i'<k$. This is equivalent to the pair $(i,k)$ lies below the ant-diagonal. We define the polynomial  
\begin{equation}\label{PXY}
P_{ik}(X,Y) = \sum_{i'\leq j\leq k} M_{ij}(X) N_{jk}(Y).
\end{equation}
 
For each  $1\leq i\leq n$, let   $D_k(X)$ stand for the lower left corner minor of order $k'$ of the matrix  $X$.  Observe that  $$D_k(X) = M_{k,k'}(X) = N_{k,k}(X).$$ 
{\bf Example}. For $n=3$ we have  
{ \small $$ X= \left(\begin{array}{ccc}
x_{11}&x_{12}&x_{13}\\
x_{21}&x_{22}&x_{23}\\
x_{31}&x_{32}&x_{33}
\end{array}\right), \qquad Y= \left(\begin{array}{ccc}
y_{11}&y_{12}&y_{13}\\
y_{21}&y_{22}&y_{23}\\
y_{31}&y_{32}&y_{33}
\end{array}\right),$$}
{ \small $$D_3(X) = x_{31},\quad D_2(X)= \left|\begin{array}{cc}
 x_{21}&x_{22}\\
 x_{31}&x_{32}
 \end{array}\right|, \quad D_1(X)=\det(X),$$
$$D_3(Y) = y_{31},\quad D_2(Y)= \left|\begin{array}{cc}
y_{21}&y_{22}\\
y_{31}&y_{32}
\end{array}\right|, \quad D_1(Y)=\det(Y),$$
$$P_{32}(X,Y) = \sum_{1\leq j\leq 2} M_{3j}(X) N_{j2}(Y)=
M_{31}(X) N_{12}(Y)+M_{32}(X) N_{22}(Y)=
x_{31} \left|\begin{array}{cc}
y_{11}&y_{12}\\
y_{31}&y_{32}
\end{array}\right| + x_{32} \left|\begin{array}{cc}
y_{21}&y_{22}\\
y_{31}&y_{32}
\end{array}\right|,$$
$$P_{33}(X,Y) = \sum_{1\leq j\leq 3} M_{3j}(X) N_{j3}(Y)= 
x_{31}y_{11}+x_{32}y_{21}+x_{33}y_{31},$$
$$P_{23}(X,Y) = \sum_{2\leq j\leq 3} M_{2j}(X) N_{j3}(Y)=
 M_{22}(X) N_{23}(Y) +   M_{23}(X) N_{33}(Y) =
\left|\begin{array}{cc}
x_{21}&x_{22}\\
x_{31}&x_{32}
\end{array}\right| y_{21} + \left|\begin{array}{cc}
x_{21}&x_{23}\\
x_{31}&x_{33}
\end{array}\right| y_{31}.$$}
\Prop\Num. The polynomials $\{P_{ik}(X,Y): ~~ i'<k\}$, and $ D_{k}(X)$,~ $ D_{k}(Y)$, where  $ 1\leq k\leq n$,  are  $U$-invariant. \\
\Proof. It is obvious that the corner  minors  $ D_{k}(X)$,~ $ D_{k}(Y)$ are   $U$-invariant.

For an arbitrary simple root  $\al=(a,a+1)$, we consider the subgroup  $s(t)=\exp(tE_{\al})$, where $t\in K$ and $E_\al=E_{a,a+1}$ is the relative matrix unit.  We denote by $\rho_\al(t)f(X)$ the action of  $s_\al(t)$ on $f(X)$. We aim to show that  $\rho_\al(t)P_{ik}(X,Y)=P_{ik}(X,Y)$, for each  $i'<k$  and    $\al=(a,a+1)$.

 We obtain 
$$\rho_\al(t) M_{ij}(X) = \left\{\begin{array}{l}M_{i,a+1}(X) +t M_{ia}(X),~~\mbox{if}~~i'< j=a+1,  \\
M_{ij} (X),~~\mbox{in~~ other~~ cases}, \end{array} \right.$$

$$\rho_\al(t) N_{jk}(Y) =\left\{\begin{array}{l} N_{ak}(Y)-tN_{a+1,k}(Y),~~\mbox{if}~~ j=a< k, \\
N_{jk}(Y), ~~\mbox{in~~ other~~ cases}.\end{array} \right.$$

Suppose that $i'\leq a<a+1\leq k$. Then
 $$\rho_\al(t)P_{ik}(X,Y) =   M_{ia}(X)  (N_{ak}(Y)-tN_{a+1,k}(Y))+ (M_{i,a+1}(X) +t M_{ia})(X) N_{a+1,k}(Y) + $$ 
$$ \sum_{\small{\begin{array}{l}i'\leq j\leq k,\\ ~j\ne a,a+1\end{array}}} M_{ij}(X) N_{jk}(Y)=$$
$$ M_{ia}(X) N_{ak}(Y)+ M_{i,a+1}(X) N_{a+1,k}(Y) + \sum_{\small{\begin{array}{l}i'\leq j\leq k,\\ ~j\ne a,a+1\end{array}}} M_{ij}(X) N_{jk}(Y)  = P_{ik}(X,Y).$$
The polynomial $P_{ik}(X,Y)$, where  $i'<k$, is $U$-invariant.

If the condition  $i'\leq a<a+1\leq k$ is not true, then all minors  $M_{ij}(X)$,~~$ N_{jk}(Y)$,  where $ i'\leq j\leq k$, are  $U$-invariant, and,  therefore, $P_{ik}(X,Y)$,~ $i'<k$, is $U$-invariant. ~~~ $\Box$

We denote 
\begin{equation}\label{PX}
 P_{ik}(X)=P_{ik}(X,X) = \sum_{i'\leq j\leq k} M_{ij}(X) N_{jk}(X).
\end{equation}
\\
\Cor\Num. The polynomials $\{P_{ik}(X):~ i'<k\}$ are $U$-invariants.
\\
\\
\Theorem\Num\label{ThUone}. The field  $K(\Mat(n))^U$ is freely generated over  $K$ by the system of polynomials  $$\{P_{ik}(X):~ i'<k\}\sqcup \{D_k(X): ~ 1\leq k\leq n\}.$$
\Proof. Let $B$ be the subgroup of upper triangular  matrices, and  $w_0$ be the element of gratest length in  the Weyl group.  

Let $\Sc$ stand for the subspace of matrices of the form  
\begin{equation}\label{SSS}
S=\small{\left(\begin{array}{cccc}
	0&0&\ldots&s_{1n}\\
	\vdots&\vdots&\ddots&\vdots\\
	0&s_{n-1,2}&\ldots&s_{n-1,n}\\
	s_{n1}&s_{n2}&\ldots&s_{nn}
	\end{array}\right). }
\end{equation}
The subset $w_0B$ is dense in  $\Sc$.
We denote by $\pi$ the restriction map   $K[\Mat(n)]^U\to K[\Sc]$.  
Since the Bruhat cell $Bw_0B$ is dense in $\GL(n)$, the subset  
$$\bigcup_{g\in U} g\Sc g^{-1}$$
is dense in  $\Mat(n)$. 
It implies that  $\pi$ is an embedding  $K[\Mat(n)]^U\hookrightarrow K[\Sc]$.
The embedding  $\pi$ extends to the embedding of fields  $$K(\Mat(n))^U\hookrightarrow K(\Sc).$$

Let $\Fc_0$ be the subfield of $K(\Mat(n))^U$ generated by  $\{D_k(X): ~ 1\leq k\leq n\}$.
 Respectively, let  $\Pc_0$ be the subfield of  $K(\Sc)$ generated by the elements $\{s_{k,k'}:1\leq k\leq n\}$ (i.e. the elements  on the anti-diagonal).
 Since $\pi(D_k(X))=\pm s_{n1}s_{n-1,2}\ldots s_{k,k'}$ for each $1\leq k\leq n$,
 the embedding $\pi$ establishes the isomorphism  $\Fc_0\to\Pc_0$.
 
For each pair $j<k$, the minor $N_{jk}(S)$  of  an arbitrary  matrix   $S$ from $\Sc$
  has the zero first row, and, therefore, it equals to zero. Then 
 $$
  \pi(P_{ik}(X)) = \pi(M_{ik}(X)N_{kk}(X)) = M_{ik}(S)N_{kk}(S). $$
  
 By definition $N_{kk}(S)=D_k(S)$. As $i'<k$,   
 $$ M_{ik}(S) = \left|\begin{array}{cccc}
 0&\ldots&0&s_{i,k}\\
 0&\ldots&s_{i+1,(i+1)'}&s_{i+1,k}\\
 \vdots&\ddots&\vdots&\vdots\\
 s_{n1}&\ldots&s_{n,(i+1)'}&s_{n,k}
 \end{array}\right|   = \pm D_{i+1}(S)s_{ik}.$$
We obtain  
  
   \begin{equation}\label{piPik}
  \pi(P_{ik}(X)) = \pm D_{i+1}(S)D_k(S)s_{ik}.
     \end{equation}
     
 The elements $\{s_{ik}: ~~i'<k\}$ freely generate the field $K(\Sc)$ over the subfield  $\Pc_0$.  The formula  (\ref{piPik}) implies that  $\pi$ 
 isomorphically maps  $K(\Mat(n))^U$ onto  $K(\Sc)$, and the system of elements  $\{P_{ik}(X):~ i'<k\}$  freely generates  the field  
  $K(\Mat(n))^U$ over $\Fc_0$. It implies the claim of the theorem. $\Box$
 
 We need the following definition.
 \\
 \Def~\Num. Let  $\{\xi_\al: \al\in \Ax\}$  and   $\{\eta_\al: \al\in \Ax\}$ be two finite systems of free generators of an extension  $F$ of the field  $K$. Let $\prec$ be a linear order  on  $\Ax$. We say that the second system of generators is obtained from the first one by a triangular transformation if  each  $\eta_\al$ can be presented in the form 
 \begin{equation}\label{etaxiorder}
 \eta_\al=\phi_\al\xi_\al+\psi_\al,
 \end{equation}
 where $\phi_\al\ne 0$ and $\phi_\al$,~ $\psi_\al$ belong to the subfield generated by  $\{\xi_\beta: ~ \beta\prec \al\}$.
 \\
 {\bf Remark}. Using the induction method, it is easy to prove that if $\{\xi_\al: \al\in \Ax\}$ is a system of free generators of a field $F$ and the other system $\{\eta_\al: \al\in \Ax\}$ is linked with the first one by formulas (\ref{etaxiorder}), then it also freely generates $F$.

 \section{$U$-invariants matrix tuples}

As in the introduction  $\Hc=\Mat(n)\oplus\ldots\oplus \Mat(n)$  is a sum of  $m$ copies of  $\Mat(n)$.  In this section, we aim to construct a system of free generators of  $K(\Hc)^U$. 
 We consider the following systems of polynomials: 
 $$\Pb_{1,\ell} =  \{P_{ik}(X_1, X_\ell):~~ 1\leq i'<k\leq n\} ~~\mbox{for~~each}~~ 2\leq\ell\leq m,$$
$$\Pb_\ell =  \{P_{ik}(X_\ell):~~ 1\leq i'<k\leq n\} ~~\mbox{for~~each}~~  1\leq \ell\leq m,$$
$$\Db_\ell =\{D_{k}(X_\ell): 1\leq k\leq n \}  ~~\mbox{for~~each}~~  1\leq \ell\leq m.$$
\Theorem\Num\label{ThU}. The field $K(\Hc)^U$ is freely generated over  $K$ by  the system of polynomials 
$$\Bb =\left(\bigcup_{\ell=2}^m \Pb_{1,\ell}\right)\cup \left(\bigcup_{\ell=1}^m \Pb_\ell\right) \cup \left(\bigcup_{\ell=1}^m \Db_\ell\right).$$
\Proof. The adjoint action  $Ad_g$,~ $g\in U$,  on $\Hc$ has the section 
$\Sc_\Hc$ that consists of matrix $m$-tuples of the form  $(S,X_2,\ldots,X_m)$, where  $S$ is a matrix from (\ref{SSS}) and $X_2,\ldots,X_m$ is an arbitrary  $(n\times n)$-matrices. 
The restriction map  $\pi:K[\Hc]^U\to K[\Sc_\Hc]$ is an embedding, and it extends to the embedding  of fields  $$\pi:K(\Hc)^U\hookrightarrow K(\Sc_\Hc).$$
Let us show that the system of polynomials  $\pi(\Bb)$ in $\Sc_\Hc$ freely generates the field  $K(\Sc_\Hc)$. 
It follows from Theorem  \ref{ThUone} that  $\pi(\Pb_1\cup\Db_1)$  freely generates  $K(\Sc)$.

We denote  $K_0=K(\Sc)$. It is sufficient to prove that for each  $2\leq \ell\leq m$  the system of polynomials  $$\pi(\Pb_{1,\ell}) \cup\Pb_\ell\cup\Db_\ell$$
freely generates the field   $K_0(X_\ell)$.  We simplify notations  $Y=X_\ell=(y_{ij})_{i,j=1}^n$. Respectively $K_0(X_\ell)=K_0(Y)$.

We denote $P_{ik}(S,Y) =\pi(P_{ik}(X,Y))$, where $S=(s_{ij})$  is a matrix of the form  (\ref{SSS}) with elements $s_{ij}=\pi(x_{ij})$.
The formula  (\ref{PXY}) takes the form  

\begin{equation}\label{PSY}
\begin{array}{l}P_{ik}(S,Y) = \sum_{i'\leq j\leq k} M_{ij}(S) N_{jk}(Y)= \\
\\
M_{i,i'}(S)N_{i',k}(Y)+ M_{i,i'+1}(S)N_{i'+1,k}(Y)+\ldots+ M_{i,k}(S)N_{k,k}(Y).\end{array}
\end{equation}

Define the linear order on the set  of pairs $\{(a,b):~ 1\leq a,b\leq m\}$ as follows 
$$(a,b)\prec (a_1,b_1),~~\mbox{if}~~ b<b_1, ~~\mbox{or}~~ b=b_1 ~~\mbox{and}~~ a>a_1.$$

Consider the subspace  $\Qc_0$ in  $K_0(Y)$ generated by $\Db_\ell=\{D_k(Y): 1\leq k\leq n\}$. The field  $\Qc_0$ is freely generated by the system  $\Db_\ell$ over $K_0$.

 Let $\Qc_1$  be an extension of the field $\Qc_0$ by the system of  polynomials  $$\Nb =\{N_{jk}(Y):~ 1\leq j < k\leq n\}.$$
 The system of polynomials $\Nb$ is algebraically independent and freely generates the field  $\Qc_1$  over  $\Qc_0$. 
 
 The elements  $\pi(\Pb_{1,\ell})$ belong to  $\Qc_1$, since all minors  $N_{jk}(Y)$,~$j\leq k$,  belong to  $\Qc_1$, and $M_{ij}(S)$,~ $i'\leq j$, belong to $K_0$.\\
 {\bf Item 1}. Let us show that  $\pi(\Pb_{1,\ell})$ freely generates the field  $\Qc_1$ over $\Qc_0$.

Consider the linear order $\prec$ on  $\Nb$.  According to the formula 
(\ref{PSY}), we get  $$ P_{ik}(S,Y) = M_{i,i'}(S)N_{i',k}(Y) + \{\mbox{terms ~~of~~ lower~~ order }\}.$$
The greatest coefficient  $M_{i,i'}(S)$ coincides with the corner minor $D_i(S)\ne 0$. 
The system of polynomials   $\pi(\Pb_{1,\ell})$ is obtained from $\Nb$ by a triangular transformation.
According to the remark at the end of previous section, the system of polynomials   $\pi(\Pb_{1,\ell})$ freely generates the field  $\Qc_1$ over  $\Qc_0$.
\\
{\bf Item 2}. Let us  show that  $\pi(\Pb_{\ell})$ freely generates  $K_0(Y)$ over $\Qc_1$. 

Easy to see that the system of matrix entries  $\Yb=\{y_{i,k}: i'<k\}$ freely generates $K_0(Y)$ over $\Qc_1$.  

The formula (\ref{PX}) implies  
\begin{equation}\label{PY}
\begin{array}{l}P_{ik}(Y) = \sum_{i'\leq j\leq k} M_{ij}(Y) N_{jk}(Y)= \\
\\
M_{i,i'}(Y)N_{i',k}(Y)+ \ldots+ M_{i,k-1}(Y)N_{k_1,k}(Y)+ M_{ik}(Y)N_{k,k}(Y).\end{array}
\end{equation}
Consider the linear order $\prec$ on $\Yb$. 
Expanding the minor  $M_{ik}(Y)$ along its first row, we get
$$M_{ik}(Y) = \pm  D_{i+1}(Y) y_{ik} + \{\mbox{expression~~over}~~ K~~\mbox{in}  ~~  y_{ab},~~ (a,b)\prec (i,k)\}.$$
Taking into account $N_{kk}(Y)=D_k(Y)$, we obtain
$$ P_{ik}(Y)= \pm D_{i+1}(Y) D_k(Y) y_{ik}  + \{\mbox{expression~~over}~~ \Qc_1~~\mbox{in}  ~~  y_{ab},~~ (a,b)\prec (i,k)\}.$$ 
It follows that the system of polynomials   $\pi(\Pb_{1,\ell})$ is obtained from $\Yb$ by a triangular transformation.
Therefore,  $\pi(\Pb_{1,\ell})$  freely generates the field  $K_0(Y)$ over $\Qc_1$. ~$\Box$

\end{document}